\documentclass[12pt,a4paper]{article}
\usepackage{amsmath,amssymb,amsbsy,amscd}
\newtheorem{thm}{Theorem}[section]
\newtheorem{lem}[thm]{Lemma}
\newtheorem{cor}[thm]{Corollary}
\newtheorem{quest}[thm]{Question}
\newtheorem{ddef}[thm]{Definition}

\newtheorem{pf}{\it Proof.}

\newcommand{\qed}{{\hfill $\Box $}}          

\newcommand{\End }{\mathop{\rm End}\nolimits}
\newcommand{\Aut }{\mathop{\rm Aut}\nolimits}
\newcommand{\BA }{\mathop{\rm BA}\nolimits}
\newcommand{\lcm}{\mathop{\rm lcm}\nolimits}
\newcommand{\zs}{\{ 0\} }
\newcommand{\sm}{\setminus}
\newcommand{\R}{{\bf R}}
\newcommand{\N}{{\bf N}}
\newcommand{\Z}{{\bf Z}}

\newcommand{\ep}{\epsilon}
\newcommand{\x}{{\bf x}}
\newcommand{\kx}{k[{\bf x}]}

\begin{document}

\title{Automorphisms of a polynomial ring 
which admit reductions of type I}

\author{Shigeru Kuroda}

\date{}

\maketitle

\begin{abstract}
Recently, Shestakov-Umirbaev solved Nagata's conjecture 
on an automorphism of a polynomial ring. 
To solve the conjecture, 
they defined notions called reductions of types I--IV 
for automorphisms of a polynomial ring. 
An automorphism admitting a reduction of type I 
was first found by Shestakov-Umirbaev. 
Using a computer, 
van~den~Essen--Makar-Limanov--Willems gave 
a family of such automorphisms. 
In this paper, 
we present a new construction of 
such automorphisms using locally nilpotent derivations. 
As a consequence, 
we discover that there exists an automorphism admitting 
a reduction of type I which satisfies some degree condition 
for each possible value. 
\end{abstract}

\section{Introduction}
\setcounter{equation}{0}

Let $k$ be a field of characteristic zero, 
and $\kx =k[x_1,\ldots ,x_n]$ the polynomial ring in $n$ variables over $k$. 
We will identify an endomorphism $F\in \End _k\kx $ 
with the $n$-tuple $(f_1,\ldots ,f_n)$ of elements of $\kx $, 
where $f_i=F(x_i)$ for each $i$. 
Then, 
$F$ is invertible if and only if $k[f_1,\ldots ,f_n]=\kx $. 
If this is the case, 
the sum $\deg F:=\sum _{i=1}^n\deg f_i$ 
of the total degrees of $f_i$'s is necessarily at least $n$. 
An automorphism $F\in \Aut _k\kx $ 
is said to be {\it affine} if $\deg F=n$, 
and {\it elementary} 
if there exist $i\in \{ 1,\ldots ,n\} $ 
and a polynomial $\phi \in \kx $ not depending on $x_i$ 
such that $f_i=x_i+\phi $ and $f_j=x_j$ for each $j\neq i$. 
We say that $F$ admits an {\it elementary reduction} 
if there exists an elementary automorphism $G$ such that 
$\deg (F\circ G)<\deg F$. 
Note that $F$ admits an elementary reduction 
if and only if there exists 
$\phi \in k[f_1,\ldots ,f_{i-1},f_{i+1},\ldots ,f_n]$ 
such that $\deg (f_i+\phi )<\deg f_i$ for some $i$. 
The subgroup $\BA _k\kx $ 
of the automorphism group $\Aut _k\kx $ 
generated by the affine automorphisms 
and the elementary automorphisms 
is called the {\it tame subgroup}, 
and each element of $\BA _k\kx $ 
is called a {\it tame automorphism}.

By Jung~\cite{Jung} and van der Kulk~\cite{Kulk}, 
it follows that $\BA _k\kx =\Aut _k\kx $ when $n=2$. 
In fact, they showed that 
each $F\in \Aut _k\kx $ for $n=2$ admits an elementary reduction 
whenever $\deg F>2$. 
Thereby, 
$$
\deg F>\deg (F\circ G_1)>\cdots >\deg (F\circ G_1\circ \cdots \circ G_r)=2 
$$
for some elementary automorphisms $G_1,\ldots ,G_r$ of $\kx $.

Now, assume that $n=3$. 
Nagata~\cite{Nagata} conjectured that the automorphism 
\begin{equation}\label{eq:Nagata}
(x_1-2(x_1x_3+x_2^2)x_2-(x_1x_3+x_2^2)^2x_3, 
x_2+(x_1x_3+x_2^2)x_3, 
x_3)
\end{equation}
of $\kx $ is not tame. 
In 2003, 
this well-known conjecture 
was finally solved in the affirmative 
by Shestakov-Umirbaev~\cite{SU1},~\cite{SU2}. 
They defined four types of reductions, 
said to be of types I, II, III and IV, 
for elements of $\Aut _k\kx $. 
Then, 
showed that each element of $\BA_k\kx $ but an affine automorphism 
admits an elementary reduction or one of these four types of reductions. 
One can easily check that Nagata's automorphism admits 
none of these reductions. 
Therefore, Nagata's conjecture is true.

\begin{ddef}[{\cite[Definition~1]{SU2}}]\label{def:typeI}\rm
Assume that $n=3$. 
We say that an automorphism 
$(f_1,f_2,f_3)$ of $\kx $ admits a {\it reduction of type I} 
if the following conditions hold, 
where we may permute the indices of $f_1$, $f_2$ and $f_3$ 
if necessary. 

(i) There exists an odd number $s\geq 3$ such that 
$\deg f_1:\deg f_2=2:s$. 

(ii) $\deg f_1<\deg f_3\leq \deg f_2$. 

(iii) $\bar{f}_3$ does not belong to $k[\bar{f}_1,\bar{f}_2]$, 
where $\bar{f}$ denotes the highest homogeneous part of $f$ 
for each $f\in \kx $. 

(iv) There exist $\alpha \in k\sm \zs $ and 
$\phi \in k[f_1,f_2-\alpha f_3]$ such that $\deg (f_3+\phi )<\deg f_3$ 
and $\deg [f_1,f_3+\phi ]<\deg f_2+\deg [f_1,f_2-\alpha f_3]$. 
Here, 
we define 
\begin{equation}\label{eq:bracket}
\deg [f,g]=\max \left\{ \deg \left(
\frac{\partial f}{\partial x_i}
\frac{\partial g}{\partial x_j}
-\frac{\partial f}{\partial x_j}
\frac{\partial g}{\partial x_i}\right) \Bigm| 1\leq i<j\leq 3\right\} +2
\end{equation}
for each $f,g\in \kx $. 
\end{ddef}

Note that $(f_1,f_2-\alpha f_3,f_3)$ 
admits an elementary reduction by (iv), 
while $(f_1,f_2,f_3)$ does not (cf.~\cite[Proposition~1]{SU2}). 
Shestakov-Umirbaev~\cite[Example~1]{SU2} 
gave the first example of a tame automorphism 
which admits a reduction of type I in case of $s=3$. 
Van~den~Essen--Makar-Limanov--Willems~\cite{EMW} 
constructed a family of such automorphisms 
when $s=3,5,7$ using a computer. 
Reductions of types II, III and IV 
are also defined theoretically~\cite{SU2}, 
but no automorphisms admitting these reductions are found. 
To study the structures of $\Aut _k\kx $ and $\BA _k\kx $, 
it is of great importance to investigate automorphisms 
admitting reductions of these four types.

The purpose of this paper is to construct 
new automorphisms of $\kx $ which admit reductions of type I 
by employing the theory of locally nilpotent derivations. 
As a consequence, 
we discover that there exists a tame automorphism 
admitting a reduction of type I such that 
$\deg f_1:\deg f_2=2:s$ for each odd number $s\geq 3$.

The author would like to thank Dr.~Hiraku Kawanoue 
for informing him of an example of polynomials 
discussed in Section~\ref{sect:remark}.

\section{Tame automorphisms admitting reductions of type I}
\setcounter{equation}{0}

Before stating our main result, 
we prove a lemma. 
In what follows, 
we assume that $n=3$.

\begin{lem}\label{lem:typeI}
Let $s\geq 3$ be an odd number, 
and $H=(h_1,h_2,h_3)$ a tame automorphism of $\kx $ 
such that 
$$
\deg h_1:\deg h_2:\deg h_3=2:s:1,\qquad 
\frac{s-1}{2}\deg h_1<\deg (ch_1^s-h_2^2)<\deg h_2
$$
for some $c\in k\sm \zs $. 
Then, 
$H'=(h_1,h_2',h_3')$ is a tame automorphism 
admitting a reduction of type I 
for which $\deg h_2'=\deg h_2$ and $\deg h_3'=\deg (ch_1^s-h_2^2)$. 
Here, 
$h_2'=h_2+h_3+ch_1^s-h_2^2$ 
and 
$h_3'=h_3+ch_1^s-h_2^2$. 
\end{lem}
\begin{pf}\rm
It is easy to check that $H'=H\circ G_1\circ G_2$, 
where $G_1$ and $G_2$ are elementary automorphisms of $\kx $ 
defined by $G_1(x_3)=x_3+cx_1^s-x_2^2$, 
$G_2(x_2)=x_2+x_3$ and 
$G_i(x_j)=x_j$ for $(i,j)\neq (1,3),(2,2)$. 
Hence, 
$H'$ is a tame automorphism of $\kx $, since so is $H$. 
By assumption, 
$\deg h_2$ is greater than $\deg h_3$ and $\deg (ch_1^s-h_2^2)$, 
while $\deg h_3$ is less than $\deg (ch_1^s-h_2^2)$. 
Hence, 
$\deg h_2'=\deg h_2$ and $\deg h_3'=\deg (ch_1^s-h_2^2)$. 
It follows that 
$$
l\deg h_1<\deg h_3'<\deg h_2'=\frac{s}{2}\deg h_1
<(l+1)\deg h_1, \text{ where }l=\frac{s-1}{2}. 
$$
This implies that 
$\bar{h}'_3$ does not belong to $k[\bar{h}_1,\bar{h}'_2]$. 
If $\alpha =1$, then $h_2'-\alpha h_3'=h_2$. 
So, 
$\phi :=-ch_1^s+h_2^2$ is contained in $k[h_1,h_2'-\alpha h_3']$. 
Then, the total degree of $h_3'+\phi =h_3$ is less than $\deg h_3'$. 
In addition, 
\begin{align*}
&\deg [h_1,h_3'+\phi ]
=\deg [h_1,h_3]
\leq \deg h_1+\deg h_3
\leq \deg h_2\\
&\quad =\deg (h_2'-\alpha h_3')<\deg (h_2'-\alpha h_3')
+\deg [h_1,\deg (h_2'-\alpha h_3')].
\end{align*}
Therefore, $(h_1,h_2',h_3')$ satisfies 
all the conditions of Definition~\ref{def:typeI}. 
\qed\end{pf}

Now, 
let $p$ and $q$ be natural numbers, 
and consider triangular derivations $D$ and $E$ of $\kx $ 
defined by 
\begin{equation}\label{eq:derivations}
\begin{gathered}
D(x_1)=x_2^{q+1}, \quad 
D(x_2)=0, \quad 
D(x_3)=(p+1)x_1^px_2^q, \\
E(x_1)=2x_3, \quad 
E(x_2)=2(p+1)x_1^p, \quad 
E(x_3)=1. 
\end{gathered}
\end{equation}
Here, 
we say that a $k$-derivation $\Delta $ of $\kx $ 
is {\it triangular} if $\Delta (x_{\sigma (i)})$ 
belongs to $k[x_{\sigma (1)},\ldots ,x_{\sigma (i-1)}]$ 
for each $i$ for some permutation $\sigma $ 
of $\{ 1,\ldots ,n\} $. 
If this is the case, $\Delta $ is {\it locally nilpotent}, 
i.e., 
$\Delta ^l(f)=0$ for sufficiently large $l$ 
for each $f\in \kx $. 
In particular, 
\begin{equation}\label{eq:polynomials}
f_i=\sum _{l=0}^{\infty }\frac{D^l(x_i)}{l!},\qquad 
g_i=\sum _{l=0}^{\infty }\frac{E^l(x_i)}{l!}(-x_3)^l
\end{equation}
are elements of $\kx $ for each $i$. 
We set $F=(f_1,f_2,f_3)$, $G=(g_1,g_2,x_3)$, 
and define $h_1=F(g_1)$, $h_2=F(g_2)$ and $h_3=f_3$. 
Namely, 
$F\circ G=(h_1,h_2,h_3)$. 
Put 
\begin{equation}\label{eq:m and c}
m=pq+p+q,\qquad  c=(-2)^{p+1}\prod _{i=1}^p\frac{i+1}{2i+1}.
\end{equation}

Here is our main result. 

\begin{thm}\label{thm:main}
Let $p$ and $q$ be natural numbers. 
Then, 
$(h_1,h_2,h_3)$ is a tame automorphism of $\kx $ for $n=3$ 
such that 
\begin{equation}\label{eq:main}
\begin{gathered}
\deg h_1=2m,\quad 
\deg h_2=(2p+1)m,\quad 
\deg h_3=m, \\ 
\deg (c^2h_1^{2p+1}-h_2^2)=2pm+p+1. 
\end{gathered}
\end{equation}
\end{thm}

Note that 
$(h_1,h_2,h_3)$ satisfies the assumptions of 
Lemma~{\rm \ref{lem:typeI}} for $s=2p+1$. 
Actually, 
$$
p\deg h_1<2pm+p+1=(2p+1)m-(p+1)q+1<\deg h_2. 
$$ 
Therefore, 
we obtain the following corollary 
to Theorem~\ref{thm:main}. 

\begin{cor}\label{cor}
There exists a tame automorphism $(h_1',h_2',h_3')$ 
of $\kx $ admitting a reduction of type I such that 
$$
\deg h_1'=2m,\quad 
\deg h_2'=(2p+1)m,\quad 
\deg h_3'=2pm+p+1 
$$
for each $p,q\in \N $, where $m=pq+p+q$. 
\end{cor}

For a triangular derivation $\Delta $, 
it is known that the exponential map $\exp \Delta :\kx \to \kx $ 
is a tame automorphism of $\kx $. 
If furthermore $\Delta (x_n)=1$, then 
$\ker \Delta =k[g_1',\ldots ,g_{n-1}']$, 
and $(g_1',\ldots ,g_{n-1}',x_n)$ is a tame automorphism of $\kx $. 
Here, we define 
\begin{equation*}
(\exp \Delta )(f)=\sum _{l=0}^{\infty }\frac{\Delta ^l(f)}{l!}, 
\qquad 
g_i'=\sum _{l=0}^{\infty }\frac{\Delta ^l(x_i)}{l!}(-x_n)^l
\end{equation*}
for each $f\in \kx $ and $i=1,\ldots ,n-1$ 
(cf.~\cite[Sections 1.3 and 6.1]{Essen}). 
Hence, 
$F$ and $G$ are tame automorphisms of $\kx $, 
and so is $F\circ G=(h_1,h_2,h_3)$. 
In addition, $E(g_i)=0$ for $i=1,2$. 
In Section~\ref{sect:proof}, 
we consider the polynomial
$$
I=x_1^{p+1}-x_2x_3. 
$$ 
Since $D(I)=0$, 
it follows that $F(I)=(\exp D)(I)=I$.

To conclude this section, we give 
explicit descriptions of $f_i$ and $g_j$ 
for $i=1,2,3$ and $j=1,2$. 
By a straightforward computation, we get 
\begin{equation}\label{eq:f1f_2g1}
f_1=x_1+x_2^{q+1},\quad f_2=x_2,\quad 
g_1=x_1-x_3^2. 
\quad 
\end{equation} 
We show that 
\begin{equation}\label{eq:f3g2}
f_3=x_3+\sum _{i=0}^p\binom{p+1}{i+1}x_1^{p-i}x_2^{(q+1)i+q},\quad 
g_2=x_2+\sum _{i=0}^pc_ix_1^{p-i}x_3^{2i+1}, 
\end{equation}
where 
\begin{equation*}
c_i=(-2)^{i+1}
\prod _{l=0}^i\frac{p-l+1}{2l+1} \quad 
\text{ for each }i. 
\end{equation*}
The first equality of (\ref{eq:f3g2}) 
is reduced to the equalities 
\begin{equation}\label{eq:f3pf}
\frac{D^{i+1}(x_3)}{(p+1)!}=\frac{x_1^{p-i}x_2^{(q+1)i+q}}{(p-i)!}
\quad \text{for }i=0,\ldots ,p. 
\end{equation}
We prove (\ref{eq:f3pf}) by induction on $i$. 
The case $i=0$ follows from the definition of $D$.  
Assume that (\ref{eq:f3pf}) is true if $i=l$ for some $0\leq l<p$. 
Then, 
\begin{align*}
&\frac{D^{l+2}(x_3)}{(p+1)!}
=D\left(\frac{D^{l+1}(x_3)}{(p+1)!}\right)=
D\left(\frac{x_1^{p-l}x_2^{(q+1)l+q}}{(p-l)!}\right) \\
&\quad =\frac{(p-l)x_1^{p-l-1}x_2^{(q+1)l+q}D(x_1)}{(p-l)!}
=\frac{x_1^{p-(l+1)}x_2^{(q+1)(l+1)+q}}{(p-(l+1))!}. 
\end{align*}
Hence, 
(\ref{eq:f3pf}) holds for $i=l+1$, and thus holds for any $0\leq i\leq p$. 
Therefore, 
we have proved the first equality of (\ref{eq:f3g2}). 
Next, 
let $g_2'$ be the right-hand side of 
the second equality of (\ref{eq:f3g2}). 
Then, $g_2-g_2'=x_2\psi $ for some $\psi \in \kx $. 
To conclude that $g_2=g_2'$, 
it suffices to show that $E(g_2'-g_2)=0$, 
since $E$ is locally nilpotent and $E(x_2)\neq 0$ by definition. 
In fact, 
for a locally nilpotent derivation $\Delta $ of $\kx $, 
the condition $\Delta (\phi \psi )=0$ implies $\psi =0$ 
for $\phi ,\psi \in \kx $ with $\Delta (\phi )\neq 0$ 
(cf.~\cite[Proposition 1.3.32]{Essen}). 
It follows that $E(g_2)=0$ as mentioned. 
Note that $c_0=-2(p+1)$, and 
$2(p-i)c_i=-(2i+3)c_{i+1}$ for $i=0,\ldots ,p$. 
Hence, we have 
\begin{align*}
E(g_2')&=E(x_2)+\sum _{i=0}^pc_i\left(
(2i+1)x_1^{p-i}x_3^{2i}E(x_3)
+(p-i)x_1^{p-i-1}x_3^{2i+1}E(x_1)
\right) \\
&=2(p+1)x_1^p+\sum _{i=0}^p\left(
(2i+1)c_ix_1^{p-i}x_3^{2i}+
2(p-i)c_ix_1^{p-i-1}x_3^{2i+2}
\right)\\
&=2(p+1)x_1^p+\sum _{i=0}^p\left(
(2i+1)c_ix_1^{p-i}x_3^{2i}-(2i+3)c_{i+1}x_1^{p-(i+1)}x_3^{2(i+1)}
\right)\\
&=2(p+1)x_1^p
+\sum _{i=0}^p(2i+1)c_ix_1^{p-i}x_3^{2i}
-\sum _{i=1}^{p+1}(2i+1)c_{i}x_1^{p-i}x_3^{2i}=0. 
\end{align*}
Thus, $E(g_2-g_2')=E(g_2)-E(g_2')=0$. 
Therefore, 
we obtain the second equality of (\ref{eq:f3g2}). 

\section{Proof of the main result}
\setcounter{equation}{0}
\label{sect:proof}
In this section, 
we prove the four equalities of (\ref{eq:main}). 

Let $f=\sum _{\alpha \in \Z ^n}c_{\alpha }\x ^{\alpha }$ 
be a Laurent polynomial in $x_1,\ldots ,x_n$ over $k$, 
where $c_{\alpha }\in k$ and 
$\x ^{\alpha }=x_1^{\alpha _1}\cdots x_n^{\alpha _n}$ 
for each $\alpha =(\alpha _1,\ldots ,\alpha _n)$. 
Then, we set 
\begin{equation*}
|f|=\{ \alpha \in \Z ^n\mid c_{\alpha }\neq 0\} . 
\end{equation*}
For $\eta \in \R ^n$, 
we define $\deg _{\eta }f$ to be the maximum 
among the inner products $\alpha \cdot \eta $ 
for $\alpha \in |f|$, and 
$$
f^{\eta }=\sum _{\alpha \in \Z ^n}c_{\alpha }'\x ^{\alpha },
\text{ \ where \ }c_{\alpha }'=\left\{ 
\begin{array}{lc}
c_{\alpha } & \text{if }\alpha \cdot \eta =\deg _{\eta }f\\
0&\text{otherwise}.
\end{array}\right.
$$
Clearly, 
$\deg _{\eta }f=\deg _{\eta }f^{\eta }$ 
for each $f\in \kx $. 
We note that $(f+g)^{\eta }$ is equal to one of 
$f^{\eta }$, $g^{\eta }$ and $f^{\eta }+g^{\eta }$ 
for each $f,g\in \kx $ 
with $|f^{\eta }|\cap |g^{\eta }|=\emptyset $. 
For a derivation $\Delta $ of $\kx $, 
we define a derivation $\Delta ^{\eta }$ of $\kx $ 
by setting 
$$
\Delta ^{\eta }(x_i)=\left\{ 
\begin{array}{cc}
(\Delta (x_i)x_i^{-1})^{\eta }x_i & 
\text{if }\deg _{\eta }(\Delta (x_i)x_i^{-1})=\deg _{\eta }\Delta \\
0&\text{otherwise}
\end{array}\right.
$$
for each $i$, 
where $\deg _{\eta }\Delta $ denotes the maximum among 
$\deg _{\eta }(\Delta (x_i)x_i^{-1})$ 
for $i=1,\ldots ,n$. 
Then, we have 
$\Delta ^{\eta }(f^{\eta })=0$ for each $f\in \ker \Delta $, 
for otherwise 
$0\neq \Delta ^{\eta }(f^{\eta })=(\Delta (f))^{\eta }$.

Now, we set $\omega _i=\deg f_i$ for $i=1,2,3$, 
and $\omega =(\omega _1,\omega _2,\omega _3)$. 
Then, 
\begin{equation*}
\omega _1=q+1,\quad 
\omega _2=1,\quad 
\omega _3=pq+p+q=m. 
\end{equation*}
For each $\alpha =(\alpha _1,\alpha _2,\alpha _3)$, 
we have 
$$
\deg F(\x ^{\alpha })=\deg f_1^{\alpha _1}f_2^{\alpha _2}f_3^{\alpha _3}
=\alpha _1\omega _1+\alpha _2\omega _2+\alpha _3\omega _3
=\deg _{\omega }\x ^{\alpha }. 
$$ 
Hence, 
$\deg F(f)\leq \deg _{\omega }f$ for each $f\in \kx $. 
The equality holds when $f^{\omega }$ is a term. 
By (\ref{eq:f1f_2g1}) and (\ref{eq:f3g2}), 
we see that $g_1^{\omega }=-x_3^2$ 
and $g_2^{\omega }=c_px_3^{2p+1}$. 
Hence, 
$\deg h_i=\deg F(g_i)=\deg _{\omega }g_i$ 
for $i=1,2$. 
Therefore, 
$\deg h_1=2m$ and $\deg h_2=(2p+1)m$. 
Besides, $\deg h_3=\deg f_3=m$. 
Thus, 
we have proved the first three equalities of (\ref{eq:main}).

Next, 
we consider the polynomial 
$P:=c^2g_1^{2p+1}-g_2^2$. 
Our goal is to establish that $\deg F(P)=2pm+p+1$, 
which immediately implies the last equality of (\ref{eq:main}). 
Write $P=P_1-P_2$, 
where 
$$
P_1=c^2g_1^{2p+1}-\phi ^2,\quad P_2=x_2^2+2\phi x_2\ \ 
\text{ and }\ \ \phi =g_2-x_2. 
$$
Set $\ep =(1,0,-2)$. 
Then, 
$g_1^{2p+1}$ and $\phi ^2$ belong to $x_3^{2(2p+1)}k[\x ^{\ep}]$, 
since $g_1$ and $\phi $ are in 
$x_3^2k[\x ^{\ep}]$ and $x_3^{2p+1}k[\x ^{\ep}]$ 
by (\ref{eq:f1f_2g1}) and (\ref{eq:f3g2}), respectively. 
Hence, 
$$
P_1^{\omega }=c'\x ^{u\ep }x_3^{2(2p+1)}=c'x_1^ux_3^{2(2p-u+1)}, 
\text{ where } u\geq 0,\ c'\in k\sm \zs . 
$$ 
We claim that $u\neq 0$. 
In fact, 
the monomial $x_3^{2(2p+1)}$ appears in 
$g_1^{2p+1}$ and $\phi ^2$ with coefficients 
$1$ and $c_p^2$, respectively. 
By definition, $c_p=c$. 
Hence, 
$x_3^{2(2p+1)}$ 
does not appear in $P_1$, so $u\neq 0$. 
On the other hand, 
$P_2^{\omega }=x_2(x_2+2\phi)^{\omega }
=2x_2\phi ^{\omega }=2cx_2x_3^{2p+1}$. 
Clearly, 
$|P_1^{\omega }|\cap |P_2^{\omega }|=\emptyset $. 
Hence, 
$P^{\omega }$ must be equal to 
$P_1^{\omega }$ or $-P_2^{\omega }$ 
or $P_1^{\omega }-P_2^{\omega }$. 
Recall that $E(g_i)=0$ for $i=1,2$. 
So, $E(P)=0$. 
This implies that 
$E^{\omega }(P^{\omega })=0$ as mentioned. 
A straightforward computation shows that 
\begin{align*}
\deg _{\omega }(E(x_1)x_1^{-1})&=\deg _{\omega }(x_1^{-1}x_3)
=-(q+1)+m=pq+q-1,\\
\deg _{\omega }(E(x_2)x_2^{-1})&=\deg _{\omega }(x_1^px_2^{-1})
=p(q+1)-1=pq+q-1,\\
\deg _{\omega }(E(x_3)x_3^{-1})&=\deg _{\omega }(x_3^{-1})
=-m=-pq-p-q<pq+q-1. 
\end{align*}
Accordingly, we get 
$E^{\omega }(x_i)=E(x_i)$ for $i=1,2$ and $E^{\omega }(x_3)=0$. 
Then, 
it follows that $E^{\omega }(P_i^{\omega })\neq 0$ for $i=1,2$. 
Therefore, 
we conclude that 
\begin{gather*}
P^{\omega }=P_1^{\omega }-P_2^{\omega }
=c'x_1^ux_3^{2(2p-u+1)}-2cx_2x_3^{2p+1}, \\
0=E^{\omega }(P^{\omega })
=2uc'x_1^{u-1}x_3^{4p-2u+3}-4c(p+1)x_1^px_3^{2p+1}, 
\end{gather*}
so $u=p+1$ and $c'=2c$. 
Consequently, we get 
$$
P^{\omega }=2cx_1^{p+1}x_3^{2p}-2cx_2x_3^{2p+1}
=2cx_3^{2p}(x_1^{p+1}-x_2x_3)
=2cx_3^{2p}I. 
$$
Hence, 
$\deg _{\omega }P=2pm+m+1$. 
Since $F(I)=I$ as mentioned, 
\begin{equation}\label{eq:degFP}
\deg F(P^{\omega })=\deg (2cg_3^{2p}I)
=2pm+p+1. 
\end{equation}
Finally, let $Q=P-P^{\omega }$. 
Since $F(P)=F(P^{\omega })+F(Q)$, 
it remains only to verify 
that $\deg F(Q)<2pm+p+1$ by (\ref{eq:degFP}). 
Note that $P$ and $Q$ belong to 
$x_2^2k[x_2^{-1}x_3^{2p+1},\x ^{\ep }]$. 
Furthermore, 
$\deg _{\omega }\x ^{\ep }=q+1-2m<0$, and 
$$
\deg _{\omega }(x_2^{-1}x_3^{2p+1})=(2p+1)m-1
=-(p+1)\deg _{\omega }\x ^{\ep }. 
$$
Hence, 
$\deg _{\omega }Q\equiv \deg _{\omega }P$ 
$\pmod{\deg _{\omega }\x ^{\ep }}$. 
Since $\deg _{\omega }Q<\deg _{\omega }P$, 
we get 
\begin{align*}
&\deg _{\omega }Q\leq \deg _{\omega }P
+\deg _{\omega }\x ^{\ep }
=2pm+m+1+q+1-2m\\
&\quad =2pm+p+1-p(q+2)+1<2pm+p+1. 
\end{align*}
Thus, 
$\deg F(Q)\leq \deg _{\omega }Q<2pm+p+1$, 
and thereby proving the last equality of (\ref{eq:main}).

\section{Remarks}\label{sect:remark}
\setcounter{equation}{0}
As far as we know, 
the answer to the following simple question is not known.

\begin{quest}\label{quest}\rm
Do there exist polynomials $f,g\in \kx $ as follows? 

(i) $k[f,g,h]=\kx $ for some $h\in \kx $. 

(ii) $\deg f:\deg g=2:3$ and $\deg (f^3+g^2)\leq \deg f$. 
\end{quest}

This question is closely related to the 
study of $\BA _k\kx $ for $n=3$. 
In fact, 
no automorphism of $\kx $ admits a reduction 
of type II or III or IV 
if the answer to Question~\ref{quest} is negative. 
The reason is as follows.

Suppose that there exists an automorphism 
of $\kx $ admitting a reduction of type II or III or IV. 
Then, 
it follows from~\cite[Definitions~2,~3~and~4]{SU2} that 
there exists an automorphism $(g_1,g_2,g_3)$ as follows: 

(1) $\deg g_1=2l$ and $\deg g_2=3l$ for some $l\in \N $. 

(2) There exists $\phi \in k[g_1,g_2]$ with $\deg \phi \leq 2l$ 
such that $\bar{\phi }$ and $\bar{g}_1$ 
are linearly independent over $k$. 

Since $\deg \phi \leq \deg g_1$ and $\deg \phi <\deg g_2$, 
the condition (2) implies that 
$\bar{\phi }\not\in k[\bar{g}_1,\bar{g}_2]$. 
Write $\phi =\sum _{i,j}c_{i,j}g_1^ig_2^j$, 
where $c_{i,j}\in k$ for each $i$ and $j$. 
Let $u_1$ and $u_2$ be the maximal numbers 
such that $c_{u_1,j'}\neq 0$ and $c_{i',u_2}\neq 0$ 
for some $j'$ and $i'$, 
and let $q_i$ and $r_i$ respectively be the quotient and residue of 
$u_i$ divided by $e_i$ for $i=1,2$. 
Here, we set
$$
e_1=\frac{\deg g_2}{\gcd (\deg g_1,\deg g_2)}=3,\quad 
e_2=\frac{\deg g_1}{\gcd (\deg g_1,\deg g_2)}=2. 
$$
Then, 
due to~\cite[Theorem 3]{SU1}, 
it follows that 
\begin{align*}
\deg \phi & \geq 
q_i(\lcm (\deg g_1,\deg g_2)-\deg g_1-\deg g_2+\deg [g_1,g_2])
+r_i\deg g_i \\
&\geq q_i(l+2)+r_i\deg g_i 
\end{align*}
for $i=1,2$. 
Since $\deg g_1=2l$ and $\deg \phi \leq 2l$ by assumption, 
$(q_1,r_1)$ must be $(0,1)$ or $(1,0)$. 
Hence, 
$u_1$ is equal to 1 or 3. 
Similarly, $(q_2,r_2)=(1,0)$ 
and so $u_2=2$, 
since $\deg g_2=3l$. 
In particular, 
$u_1\leq 3$ and $u_2=2$. 
The polynomials 
$g_1^ig_2^j$ for $i=0,1,2,3$ and $j=0,1,2$ 
but $g_1^3$ and $g_2^2$ have distinct total degrees. 
This implies that $c_{i,j}=0$ 
for each $(i,j)$ with $2i+3j>6$, 
while $c_{3,0}\neq 0$ and $c_{0,2}\neq 0$. 
In fact, otherwise 
$\bar{\phi }=c_{i,j}\bar{g}_1^i\bar{g}_2^j$ 
for some $(i,j)$. 
This contradicts that $\bar{\phi }\not\in k[\bar{g}_1,\bar{g}_2]$. 
Hence, 
$c_{3,0}\neq 0$, $c_{0,2}\neq 0$ and 
\begin{equation}\label{eq:ell0}
\phi =c_{3,0}g_1^3+c_{0,2}g_2^2
+c_{1,1}g_1g_2+c_{2,0}g_1^2+c_{0,1}g_2
+c_{1,0}g_1+c_{0,0}. 
\end{equation}
Without loss of generality, 
we may assume that $c_{0,2}=1$. 
Then, (\ref{eq:ell0}) is expressed as 
\begin{equation}\label{eq:ell1}
\phi =c_{3,0}\hat{f}^3+\hat{g}^2+b\hat{f}+c,
\text{ where } \hat{f}=g_1+a,\ \hat{g}=g_2+\frac{c_{1,1}}{2}g_1+\frac{c_{0,1}}{2}\end{equation}
and $a,b,c\in k$. 
Indeed, 
$\phi =c_{3,0}g_1^3+g^2+c_{2,0}'g_1^2+c_{1,0}'g_1+c_{0,0}'$ 
for some $c_{2,0}',c_{1,0}',c_{0,0}'\in k$. 
Then, we get (\ref{eq:ell1}) for $a=c_{0,2}'/3$. 
Finally, put $f=c_{3,0}\hat{f}$ and $g=c_{3,0}\hat{g}$. 
Clearly, 
$\deg f=\deg g_1=2l$, $\deg g=\deg g_2=3l$, 
and $k[f,g,g_3]=\kx $. 
Moreover, 
$$
\deg (f^3+g^2)=\deg c_{3,0}^2(c_{3,0}\hat{f}^3+\hat{g}^2)\leq 2l=\deg f
$$
by (\ref{eq:ell1}), 
since the total degrees of $\deg \phi$ and $b\hat{f}+c$ 
are at most $2l$. 
Therefore, $f$ and $g$ satisfy the 
conditions of Question~\ref{quest}.

It is worthwhile to mention that, 
if there exists a tame automorphism $(h_1,h_2,h_3)$ 
with $\deg h_1:\deg h_2:\deg h_3=2:3:1$ 
and $\deg (ch_1^3-h_2^2)\leq \deg h_1$ for some $c\in k\sm \zs $, 
then we can construct a tame automorphism 
admitting a reduction of type II or III. 
On the other hand, 
(\ref{eq:main}) gives that 
$\deg h_1=2m$, $\deg h_2=3m$, $\deg h_3=m$ 
and $\deg (c^2h_1^3-h_2^2)=2m+2$ if $p=1$. 
In this case, we have 
$$
\frac{\deg h_1}{\deg (c^2h_1^3-h_2^2)}=\frac{2m}{2m+2}\to 1
\qquad (q\to \infty ), 
$$
although $\deg (c^2h_1^3-h_2^2)>\deg h_1$.

Assume that $f,g\in \kx $ are algebraically independent over $k$ 
for which $\deg f:\deg g=r:s$, 
where $r,s\in \N $ with $2\leq r<s$ and $\gcd (r,s)=1$. 
Then, 
it easily follows from \cite[Theorem 3]{SU1} that 
$$
\deg (f^s+g^r)>\left\{ 
\begin{array}{cl}
\deg g & \text{if }r\geq 3 \\
\deg f & \text{if }r=2\text{ and }r\geq 5. 
\end{array}
\right.
$$
Hence, 
$\deg (f^s+g^r)\leq \deg f$ is possible only if $(r,s)=(2,3)$. 
We define $f,g\in \kx $ by
\begin{equation}\label{eq:kawanoue}
f=-x_1^{4l}x_2^{2(2m-1)}-2x_1^lx_2^m,\quad 
g=x_1^{6l}x_2^{3(2m-1)}+3x_1^{3l}x_2^{3m-1}+\frac{3}{2}x_2
\end{equation}
where $l,m\in \N $. 
Then, 
it is easy to see that 
$\deg f:\deg g=2:3$, 
$f$ and $g$ are algebraically independent over $k$, 
and 
$$
f^3+g^2=x_1^{3l}x_2^{3m}+\frac{9}{4}x_2^2. 
$$
In particular, 
$\deg (f^3+g^2)=\deg f$ if $l=m=1$, 
and $\deg (f^3+g^2)<\deg f$ otherwise. 
If $k$ is of characteristic $r>0$, then 
$f=x_1^{rl}$ and $g=x_2+x_1^{sl}$ satisfy 
$f^s-g^r=-x_2^r$ for any $l,s\in \N $. 
Therefore, 
$\deg (f^s-g^r)\leq \deg f$ in this case. 

Note: 
Instead of Question~\ref{quest}, 
the author first asked a question whether there exist 
$f,g\in \kx $ with $\deg f:\deg g=2:3$ and 
$\deg (f^3+g^2)\leq \deg f$ 
which are algebraically independent over $k$. 
In answer to the question, 
Dr.~Hi\nolinebreak raku Kawanoue informed him of an example 
with $\deg (f^3+g^2)=\deg f$. 
The example (\ref{eq:kawanoue}) is a modification of 
Kawanoue's example.

\begin{flushleft}
Department of Mathematics and Information Sciences\\ 
Tokyo Metropolitan University\\
1-1  Minami-Ohsawa, Hachioji \\
Tokyo 192-0397, Japan\\
E-mail: {\tt kuroda@tmu.ac.jp}
\end{flushleft}

\end{document}